# The Humbert-Bessel functions, Stirling numbers and probability distributions in coincidence problems

D. Babusci[1], G. Dattoli[2], E. Di Palma[2], E. N. Petropoulou[3]

[1]*INFN- Laboratori Nazionali di Frascati, via E. Fermi 40-00044 Frascati, Italy*

[2]*ENEA-Centro Ricerche Frascati, via E. Fermi 45- 00044 Frascati, Italy*

[3]*Department of Engineering Sciences, University of Patras, 26500 Patras, Greece*

### Abstract

The Humbert-Bessel are multi-index functions with various applications in electromagnetism. New families of functions sharing some similarities with Bessel functions are often introduced in the mathematical literature, but at a closer analysis they are not new, in the strict sense of the word, and are shown to be expressible in terms of already discussed forms. This is indeed the case of the re-modified Bessel functions, whose properties have been analyzed within the context of coincidence problems in probability theory. In this paper we show that these functions are particular cases of the Humbert-Bessel ones.

## 1 Introduction

The two index Humbert-Bessel functions have been introduced in Mathematics during the first half of the last century [1]. They can be defined through the generating function [2]

$$\sum_{m_1,m_2=-\infty}^{+\infty} u^{m_1} v^{m_2} I_{m_1,m_2}(x) = e^{u+v+\frac{x}{uv}} \qquad (1)$$

and are expressed by the series

$$I_{m_1,m_2}(x) = \sum_{r=0}^{\infty} \frac{x^r}{r!(m_1+r)!(m_2+r)!} \qquad (2).$$

The relevant recurrence relations can be easily derived from the previous relations and read

$$\partial_x I_{m_1,m_2}(x) = I_{m_1+1,m_2+1}(x),$$
$$m_1 I_{m_1,m_2}(x) = I_{m_1-1,m_2}(x) - x I_{m_1+1,m_2+1}(x),$$
$$m_2 I_{m_1,m_2}(x) = I_{m_1,m_2-1}(x) - x I_{m_1+1,m_2+1}(x) \quad (3)$$

By combining the second and first recurrences, we find

$$(m_k + x\partial_x) I_{m_1,m_2}(x) = I_{m_1-\delta_{1,k}, m_2-\delta_{2,k}}(x), k=1,2 \quad (4)$$

and, accordingly, we define the following ladder operators

$$\hat{E}_k^- = \hat{N}_k + x\partial_x, k=1,2$$
$$\hat{E}_{1,2}^+ = \hat{E}^+ = \partial_x \quad (5)$$

with $\hat{N}_k$ being a number operator such that

$$\hat{N}_k I_{m_1,m_2}(x) = m_k I_{m_1,m_2}(x) \quad (6),$$

thus obtaining

$$\hat{E}_k^- I_{m_1,m_2}(x) = I_{m_1-\delta_{1,k}, m_2-\delta_{2,k}}(x),$$
$$k=1,2 \quad (7).$$

The differential equation satisfied by these functions can thereby written as (we omit the argument of the function for conciseness)

$$\hat{E}_1^- \hat{E}_2^- \hat{E}^+ I_{m_1,m_2} = I_{m_1,m_2},$$
$$\hat{E}^+ \hat{E}_1^- \hat{E}_2^- I_{m_1,m_2} = I_{m_1,m_2} \quad (8),$$

which in explicit form writes

$$(m_1 + 1 + x\partial_x)(m_2 + 1 + x\partial_x)\partial_x I_{m_1,m_2} = I_{m_1,m_2},$$
$$\partial_x (m_1 + x\partial_x)(m_2 + x\partial_x) I_{m_1,m_2} = I_{m_1,m_2} \quad (9),$$

i.e., they satisfy the third order differential equation with non-constant coefficients ($I_{m_1,m_2} = z$)

$$x^2 z''' + (m_1 + m_2 + 3) x z'' + (m_1 + 1)(m_2 + 1) z' - z = 0 \quad (10).$$

The Humbert-Bessel functions can be generalized to more general forms including many (> 2) indices [2] and it will be shown that they are the natural framework to place the theory of the re-modified Bessel functions [3].

## 2 Multi-index and remodified Bessel functions

The functions

$$I_0(n, x) = \sum_{r=0}^{\infty} \left(\frac{x}{n}\right)^{nr} \frac{1}{(r!)^n} \qquad (11)$$

have been recently introduced [3] in probabilistic theories associated with coincidence and have been called re-modified Bessel. The reason of such a definition is simply due to the fact that for n=2 they reduce the 0-th order modified Bessel functions of first kind.

The above functions have very interesting properties, useful for the specific applications of ref. [3], and therefore their study is worth to be pursued. We will prove that their understanding in terms of already defined functions of Bessel type may significantly simplify such a task.

The case n=3 can be expressed in terms of Humbert-Bessel functions as

$$I_0(3, x) = I_{0,0}\left(\left(\frac{x}{3}\right)^3\right) \qquad (12)$$

Therefore the relevant differential equation can be derived from eq. (9) as

$$\partial_\xi (\xi \partial_\xi)^2 y = y, \quad y = I_0(3, x)$$

$$\xi = \left(\frac{x}{3}\right)^3, \qquad (13)$$

which, in terms of the x-variable, reads (see also ref. [3])

$$x^2 y''' + 3 x y'' + y' - x^2 y = 0 \qquad (14).$$

A generalization of the Humbert functions is provided by the generating function

$$\sum_{\{m\}=-\infty}^{+\infty} \{u\}^{\{m\}} I_{\{m\}}(x) = e^{\sum_{j=1}^{p} u_j + \frac{x}{\prod_{k=1}^{p} u_k}},$$

$$\{m\} = m_1, \ldots, m_p, \qquad (15)$$

$$\{u\} = u_1, \ldots, u_p,$$

$$I_{\{m\}}(x) = I_{m_1, \ldots, m_p}(x)$$

and are defined by the series

$$I_{\{m\}}(x) = \sum_{r=0}^{\infty} \frac{x^r}{r! \prod_{k=1}^{p}(m_k + r)!} \qquad (16).$$

The ladder operators are easily introduced in this case too, and the relevant differential equation reads

$$\partial_x \prod_{k=1}^{p}(m_k + x\partial_x) z = z, \qquad (17)$$
$$z = I_{\{m\}}(x)$$

In the case in which $\{m\} = 0$, the previous identity reduces to

$$\partial_x [x\partial_x]^p z = z \qquad (18).$$

Furthermore the use of the operational identity [4]

$$[x\partial_x]^p = \sum_{r=0}^{p} S_2(r,p) x^r \partial_x^r \qquad (19)$$

with $S_2(r,p)$ being the Stirling numbers of the second kind, yields

$$\partial_x \left( \sum_{r=0}^{p} S_2(r,p) x^r \partial_x^r \right) z = z \qquad (20).$$

It is now evident that

$$I_0(n,x) = I_{\{0\}}\left(\left(\frac{x}{n}\right)^n\right) \qquad (21),$$

for $p = n - 1$ and, since

$$\xi \partial_\xi = \frac{x}{n} \partial_x$$
$$\xi = \left(\frac{x}{n}\right)^n \qquad (22),$$

the differential equation of the re-modified Bessel reads

$$\partial_x \left( \sum_{r=0}^{n-1} S_2(r, n-1) x^r \partial_x^r \right) y = x^{n-1} y, \qquad (23).$$
$$y = I_0(n,x)$$

The q-th order re-modified Bessel are defined as [3]

$$I_q(n, x) = \left(\frac{x}{n}\right)^q \sum_{r=0}^{\infty} \left(\frac{x}{n}\right)^{nr} \frac{1}{(r!)^n (r+1)^q}, \qquad (24).$$

$q \leq n-1$

According to the previous definition we get[1]

$$I_q(n, x) = \left(\frac{x}{n}\right)^q I_{\{m\}}\left(\left(\frac{x}{n}\right)^n\right), \qquad (25)$$

$\{m\} = m_1 = 1, ..., m_q = 1, m_{q+1} = 0, ..., m_{n-1} = 0$

and the relevant equation can accordingly be written as

$$\partial_x \sum_{r=0}^{q} \binom{q}{r} \frac{n^q}{n^r} (x\partial_x)^{r+n-q-1} \left(x^{-q} w\right) = x^{n-(q+1)} w$$

$$w = I_{\{m\}}\left(\left(\frac{x}{n}\right)^n\right) \cdot \left(\frac{x}{n}\right)^q = I_q(n, x) \quad \text{with} \quad q < n-1$$

(26)

which, on account of eq. (19), can also be written as

$$\partial_x \sum_{r=0}^{q} \binom{q}{r} \frac{n^q}{n^r} \sum_{s=0}^{r+n-q-1} S_2(s, r+n-q-1) x^s \partial_x^s \left(x^{-q} w\right) = x^{n-(q+1)} w \qquad (27).$$

The formula (26) is derived directly from the relation (17) without using other recursive relations; the result is that the ODE equation for (25) is one degree below of that reported in ref[3] for q=1.

The operational methods and the use of the properties of the Humbert-Bessel functions greatly simplify the study of the Bessel functions of ref. [3]. Further comments will be presented in the forthcoming section.

### 3 Final Comments

In a recent series of papers [5] the theory of Bessel functions has been developed using a symbolic method which revealed particularly useful, for example, in dealing with the evaluation of the associated integrals.

---

[1] It should be noted that in the case of n=2, we are left with q=1 and with

$$I_1(2, x) = \sum_{r=0}^{\infty} \left(\frac{x}{2}\right)^{2r+1} \frac{1}{k!(k+1)!} = I_1(x)$$

$$I_n(x) = \sum_{r=0}^{\infty} \left(\frac{x}{2}\right)^{2r+n} \frac{1}{k!(k+1)!}$$

The same method can be extended to the Humbert-Bessel functions, which, within such a context, can be defined as

$$I_{m_1,m_2}(x) = \hat{c}_1^{m_1} \hat{c}_2^{m_2} e^{\hat{c}_1 \hat{c}_2 x} \varphi_1(0) \varphi_2(0) \qquad (28)$$

where $\hat{c}_k$ are operators, which, acting on the "vacuum" $\varphi_k(0)$, yield

$$\hat{c}_k^{v_k} \varphi_k(0) = \frac{1}{\Gamma(m_k + 1)}, \qquad (29).$$
$$v_k \in R$$

Just to give a flavor of the usefulness of the above redefinition we note that integrals of the type

$$L(\beta) = \int_0^\infty I_{0,0}(-x) e^{-\beta x} dx, \qquad (30)$$

can easily be derived. By treating indeed the operators $\hat{c}$ as ordinary numbers we find[2]

$$L(\beta) = \frac{1}{\beta + \hat{c}_1 \hat{c}_2} \varphi_1(0)\varphi_2(0) =$$
$$= \frac{1}{\beta} \sum_{r=0}^\infty (-1)^r \left(\frac{\hat{c}_1 \hat{c}_2}{\beta}\right)^r \varphi_1(0)\varphi_2(0) = \qquad (31)$$
$$= \frac{1}{\beta} J_0\left(\frac{2}{\sqrt{\beta}}\right)$$

where $J_0(x) = \sum_{r=0}^\infty \frac{(-1)^r x^{2r}}{2^{2r} (r!)^2}$ is the 0-th order cylindrical Bessel function of the first kind.

As a further example, it's easy to show that

$$\int_{-\infty}^\infty I_{0,0}(x) e^{-\beta x^2} dx = \sqrt{\frac{\pi}{\beta}} I_{0,0}\left(\frac{1}{4\beta} \mid 2\right),$$
$$I_{m_1,m_2}(x \mid k) = \sum_{r=0}^\infty \frac{x^r}{r! \Gamma(k r + 1 + m_1) \Gamma(k r + 1 + m_2)} \qquad (32)$$

According to the above definition we also get

---
[2] The series are summed regardless any convergence condition and they are assumed to be formally true

$$I_0(3,x) = e^{\hat{C}\left(\frac{x}{3}\right)^3} \Phi(0) \qquad (33)$$

with

$$\hat{C}^m \Phi(0) = \frac{1}{[\Gamma(m+1)]^2} \qquad (34).$$

We have, accordingly, reduced the re-modified Bessel to a third order exponential function. In this specific case the use of the following identity from the Airy transform theory [6]

$$e^{\lambda x^3} = \int_{-\infty}^{+\infty} dt\, e^{\sqrt[3]{\lambda} x t} Ai(t)$$

$$Ai(t) = \frac{1}{2\pi} \int_{-\infty}^{+\infty} \exp\left[\frac{i}{3}\xi^3 + it\xi\right] d\xi \qquad (35)$$

yields

$$I_0(3,x) = \int_0^\infty dt\, e^{\sqrt[3]{\hat{C}}\left(\frac{xt}{3}\right)} Ai(t) \Phi(0) = \int_0^\infty dt\, I_{0,0}\left(\frac{xt}{3} \mid \frac{1}{3}\right) Ai(t) \qquad (36).$$

The use of other generalized transform can be exploited to provide an alternative formulation of the Bessel re-modified functions. This aspect of the problem will be discussed in a forthcoming publication.

### References

[1] P. Humbert, Sur les fonctions du troisiemme order, C. R. A. S. 190, 159-160 (1930)

[2] G. Dattoli, S. Lorenzutta, G. Maino, G.K. Voykov and C. Chiccoli, Theory of two-index functions and applications to physical problems, J. Math. Phys. 35, 3636 (1994)

[3] M. Griffiths, Remodified Bessel Functions via Coincidences and Near Coincidences, J. Integer Seq. , 14 (2011), Art. 11.7.1

[4] W. Lang, On generalizations of the Stirling Number Triangles, J. Integer Seq.,3 (2000), Art. 00.2.4 ,pp. 1-11

[5] K. Gorska, D. Babusci, G. Dattoli, G. H. E. Duchamp and K. Penson, The Ramanujan Master Theorem and its implications for Special Functions, arXiv:1104.3406v1 [math-ph] 2011 (submitted for publication).

[6] O. Vallée and M. Soares, Airy Functions and application to Physics, World Scientific, London (2004).